 \documentclass[12pt]{article}
 \usepackage{amsmath,amsfonts,theorem}
 \numberwithin{equation}{section}

 \newtheorem{prop}{Proposition}[section]

 \newtheorem{cor}{Corollary}[section]

 \newtheorem{thm}{Theorem}[section]

 \theorembodyfont{\rmfamily}
 \newtheorem{dfn}{Definition}[section]

 \newtheorem{rmk}{Remark}

 \newcommand{\1}{^{-1}}
 \newcommand{\dual}{^\vee}
 \newcommand{\iso}{\cong}
 \newcommand{\Span}[1]{\left< #1 \right>}
 \newcommand{\half}{\frac12}
 \newcommand{\dd}{\mathrm{d}}
 \newcommand{\ev}{\mathrm{ev}}
 \newcommand{\odd}{\mathrm{odd}}
 \newcommand{\id}{\operatorname{id}}
 \newcommand{\codim}{\operatorname{codim}}
 \newcommand{\CS}{\operatorname{CS}}
 \newcommand{\Ab}{_{\mathrm{Ab}}}
 \newcommand{\nonAb}{_{\mathrm{non-Ab}}}
 \newcommand{\BS}{\mathrm{BS}}
 \newcommand{\BPU}{\mathrm{BPU}}
 \newcommand{\hWBS}{\mathrm{hWBS}}
 \newcommand{\hF}{\mathrm{hF}}
 \newcommand{\kBS}{k\mathrm{\text{-}BS}}
 \newcommand{\KS}{\mathrm{KS}}
 \newcommand{\LC}{\mathrm{LC}}
 \newcommand{\semis}{^{\mathrm{ss}}}

 \newcommand{\irr}{^{\mathrm{irr}}}
 \newcommand{\red}{\mathrm{red}}
 \newcommand{\triv}{\mathrm{triv}}
 \newcommand{\rest}[1]{{}_{{\textstyle{|}}#1}} % restriction of map

 \newcommand{\Dbar}{\overline{D}}
 \newcommand{\phibar}{\overline{\fie}}
 \newcommand{\Debar}{\overline{\De}}
 \newcommand{\hFbar}{\overline{\hF}}
 \newcommand{\wave}{\widetilde}
 \newcommand{\tensor}{\otimes}
 \newcommand{\wBS}{\wave{\BS}}
 \newcommand{\wG}{\wave{G}}
 \newcommand{\wGa}{\wave{\Ga}}
 \newcommand{\ws}{\wave{s}}
 \newcommand{\wsL}{\wave{\sL}}
 \newcommand{\wMod}{\wave{\Mod}}
 \newcommand{\wSi}{\wave{\Si}}
 \newcommand{\wW}{\wave{W}}

 \newcommand{\cosum}{\mathrel{\#}}
 \newcommand{\sI}{\mathcal I}
 \newcommand{\sH}{\mathcal H}
 \newcommand{\sL}{\mathcal L}
 \newcommand{\sM}{\mathcal M}
 \newcommand{\Oh}{\mathcal O}

 \newcommand{\al}{\alpha}
 
 \newcommand{\ga}{\gamma}
 \newcommand{\de}{\delta}
 \newcommand{\fie}{\varphi}
 
 \newcommand{\om}{\omega}
 \newcommand{\De}{\Delta}
 \newcommand{\Ga}{\Gamma}
 \newcommand{\La}{\Lambda}
 \newcommand{\Om}{\Omega}
 \newcommand{\Si}{\Sigma}

 \newcommand{\PP}{\mathbb P}
 \newcommand{\C}{\mathbb C}

 \newcommand{\R}{\mathbb R}
 \newcommand{\Z}{\mathbb Z}

 \newcommand{\res}{\operatorname{res}}
 \newcommand{\rk}{\operatorname{rank}}
 \newcommand{\corank}{\operatorname{corank}}
 \newcommand{\Aut}{\operatorname{Hom}}
 
 \newcommand{\tr}{\operatorname{tr}}
 \newcommand{\End}{\operatorname{End}}
 \newcommand{\LaGr}{\operatorname{\La_{\uparrow}\!}}
 
 \newcommand{\Hom}{\operatorname{Hom}}
 
 \newcommand{\Mod}{\operatorname{Mod}}
 \newcommand{\Mon}{\operatorname{Mon}}
 \newcommand{\Op}{\operatorname{Op}}

 \newcommand{\Sing}{\operatorname{Sing}}
 \newcommand{\Vol}{\operatorname{Vol}}

 \newcommand{\SO}{\operatorname{SO}}
 \newcommand{\Spin}{\operatorname{Spin}}
 \newcommand{\U}{\operatorname U} % unitary group
 \newcommand{\PU}{\operatorname{PU}}
 \renewcommand{\Sp}{\operatorname{Sp}}
 \newcommand{\SU}{\operatorname{SU}}
 \newcommand{\fsl}{\operatorname{\mathfrak{sl}}}

\begin{document}

 \title{On Bohr--Sommerfeld bases}
 \markright{\hfill On Bohr--Sommerfeld bases \quad}

 \author{Andrei Tyurin}
 \date{September 3, 1999}
 \maketitle

 \begin{abstract}
This paper combines algebraic and Lagrangian geo\-metry to construct a
special basis in every space of conformal blocks, the Bohr--Sommerfeld
(BS) basis. We use the method of Borthwick--Paul--Uribe \cite{BPU}, whereby
every vector of a BS basis is defined by some half-weighted Legendrian
distribution coming from a Bohr--Sommerfeld fibre of a real polarization
of the underlying symplectic manifold. The advantage of BS bases (compared
to bases of theta functions in \cite{T1}) is that we can use information
from the skillful analysis of the asymptotics of quantum states. This gives
that Bohr--Sommerfeld bases are unitary quasi-classically. Thus we can
apply these bases to compare the Hitchin connection {\cite{H}} with the KZ
connection defined by the monodromy of the Knizhnik--Zamolodchikov
equation in combinatorial theory (see, for example, Kohno \cite{K1} and
\cite{K2}).
 \end{abstract}

 \section{Degree 0 cycles}
 \markright{\hfill On Bohr--Sommerfeld bases \quad}

Let $X$ be a K\"ahler manifold with K\"ahler form $\om$ and $L$ an
algebraic geometric polarization with $c_1(L)=[\om]\in H^2(X,\Z)$. Suppose
in addition that the canonical class is even:
 \[
K_X = 0 \mod2,
 \]
and fix a {\em metaplectic structure} on $X$, that is, a line bundle
$L_{K/2}$ such that $L_{K/2}^{\tensor2}=L_{K_X}$. Then $L^k\tensor L_{K/2}$
is a holomorphic line bundle on $X$ for any $k\in\Z^+$, here called the
{\em level}. We get spaces
 \begin{equation}
 \sH_I^k=H^0(X,L^k\tensor L_{K/2}) \quad\text{for $k\ge0$}
 \label{eq1.1}
 \end{equation}
 of holomorphic sections of $L^k\tensor L_{K/2}$; here $I$ denotes the
complex structure of $X$. Thus if $\sM$ is a family of polarized complex
structures on the underlying smooth compact manifold $X$, the spaces
(\ref{eq1.1}) define holomorphic vector bundles
 \begin{equation}
 \sH^k\to \sM,
 \label{eq1.2}
 \end{equation}
 with fibres (\ref{eq1.1}) for $k\gg0$.

For the applications we have in mind, it is sufficient to work under the
following restriction: there exists an integer $d$ such that 
 \begin{equation}
 L_{K/2}=L^d. 
 \label{eq1.3}
 \end{equation}
Thus twisting by $L_{K/2}$ just shifts the level $k$ to $k+d$. We
use the hyperplane section class
 \[
 c_1(L)=[\om]\in H^2(X,\Z),
 \]
to define the {\em degree} of a cycle $C\subset X$: it is the integer
 \[
 \deg C=[C]\cdot[\om]^{(\dim C)/2}=\int_C(\om\rest C)^{(\dim C)/2}.
 \]
For any holomorphic (effective algebraic) subcycle $C$
 \[
 \dim C>0\implies \deg C > 0.
 \]
 More precisely, there is an exact sequence
 \[
 0\to \sI_C\to \Oh_X\to \Oh_C\to 0,
 \]
and restricting global sections (constants) to $C$ defines a distinguished
line in the space of sections:
 \begin{equation}
\C=H^0(\Oh_X)\to H^0(\Oh_C).
 \notag
 \label{eq1.4}
 \end{equation}

Our line bundle $L^{k+d}$ defines a rational (that is, meromorphic) map
$X\to\PP(\sH_I^{k+d})\dual$, sending $x\in X$ to the hyperplane
$H^0(\sI_x\tensor L^{k+d})\subset H^0(L^{k+d})$ of sections vanishing at
$x$. For an open set $U$, a trivialization of $L\rest{U}$ is given by a
section $s_U$ that is everywhere nonvanishing on $U$, and $x\in U$ defines
a covector of $H^0(L^k)$ of evaluation at $x$, using $s_U$ to identify the
fibre at $x$ with $\C$. More technically, there is an exact sequence
 \begin{equation}
0\to \sI_x\tensor L^{k+d}\to L^{k+d}\to \Oh_x\to 0
\label{eq1.5}
 \end{equation}
where the epimorphism is the restriction homomorphism. Part of its
cohomo\-logy sequence
 \begin{equation}
 H^0(\sI_x\tensor L^{k+d})\to H^0(L^{k+d})\to \C\to H^1(\sI_x\tensor L^{k+d})
 \label{eq1.6}
 \end{equation}
shows that $H^0(\sI_x\tensor L^{k+d})$ is indeed a hyperplane if
$H^1(\sI_x\tensor L^{k+d})=0$. Now by Serre's classical Theorems~A and~B,
$H^1(\sI_x\tensor L^{k+d})=0$ for $k\gg0$, and we have a map
 \begin{equation}
 \PP \fie_k\colon X\to \PP H^0(L^{k+d}),
 \label{eq1.7}
 \end{equation}
which is an embedding for $k\gg0$. This standard construction of algebraic
geometry reduces the study of $X$ to projective geometry.

 The holomorphic structure of $L$ admits a Hermitian connection $a_L$,
defined by the complex structure, with curvature form $2\pi i\om$. On the
other hand, our K\"ahler metric defines a Levi-Civita Hermitian connection
$a_{\LC}$ on $L_{K/2}$. We suppose also that
 \begin{equation}
 a_{L^d}=a_{\LC}
 \label{eq:aLC}
 \end{equation}
(see (1.3)).

We can now define a Hermitian form on $H^0(L^{k+d})$ in two steps: first,
every section $\ws$ of $L^{k+d}$ is locally of the form $s\cdot\hF$, where
$s$ is a section of $L^k$ and $\hF$ a section of $L_{K/2}$ (a {\em
half-form}). For two such sections $\ws_1=s_1\cdot\hF_1$ and
$\ws_2=s_2\cdot\hF_2$, set 
 \begin{equation}
 \Span{\ws_1,\ws_2}=\int_X(s_1,s_2)\cdot (\hF_1,\hF_2).
 \label{eq1.8}
 \end{equation}
 We get an identification of vectors and covectors:
 \[
 H^0(L^{k+d})=\overline{H^0(L^{k+d})}^*.
 \]
 In particular, a trivialization over an open $U$ sends
 \[
 \fie_U\colon U\to H^0(L^{k+d})
 \]
 and the projectivization of this is just the complex conjugate of (1.6):
 \begin{equation}
 \PP \fie_U\colon U\to \PP H^0(L^{k+d}).
 \label{eq1.9}
 \end{equation}
 In the set-up of complex quantization, vectors in $H^0(L^{k+d})$ are
called {\em states}, and vectors $\phibar(x)$ for $x\in X$ are called {\em
coherent states} (recall that states are {\em distributions}, not
functions).

From now on we can forget about what $K/2$ means geometrically, and
consider a twisting by the metaplectic structure as a shift of level.

Now inverting the usual way of thinking, we ask whether there are
submanifolds of $X$ of degree 0 that define hyper\-planes in the spaces of
sections $H^0(X,L^k)$, and what kind of submanifolds these are. We
strengthen the condition $\deg C=0$ to $\om\rest{C}=0$; in other words,
every such submanifold $\sL$ must be isotropic with respect to the
K\"ahler form $\om$. Thus a maximal dimensional submanifold must be
Lagrangian. Just as an algebraic subvariety may be singular, we do not
need to restrict ourselves to Lagrangian submanifolds: in what follows we
consider Lagrangian cycles a priori admitting singularities. The main
property of any such cycle is
 \begin{quote}
 {\em it can't be contained in any proper algebraic subvariety}
 \end{quote}
(in particular, in a divisor). Thus any holomorphic object is uniquely
determined by its restriction to a Lagrangian cycle $\sL$. Thus
restrictions to $\sL$ can serve as boundary conditions for holomorphic
sections of line bundles with curvature proportional to $\om$.

 \begin{rmk}
 Geometrically, if we consider Lagrangian cycles as supports of boundary
conditions for holomorphic objects, they have the minimal possible
dimension. Usual boundary conditions deal with boundaries of complex
domains of real codimension 1. Thus it is only for Riemann surfaces that
Lagrangian boundary conditions coincide with the usual boundary
conditions. In this case, in the modern theory of integrable systems the
restriction of holomorphic objects to a small circle around a point
reduces many analytical problems to algebraic geometry of curves (see for
example the survey \cite{DKN}). Thus we should add the role of Lagrangian
submanifolds as boundary conditions for holomorphic objects to Alan
Weinstein's proclamation \cite{Wei1}, p.~5. It seems reasonable to expect
that restrictions to Lagrangian submanifolds give a higher dimensional
generalization of the modern version of the theory of integrable systems. 
 \end{rmk}

If we forget for a minute the complex structure $I$ on $X$, the
polarization $L$ gives us a quadruple
 \begin{equation}
 (X,\om,L,a_L),
 \end{equation}
 where $\om$ is the K\"ahler form and $a_L$ a Hermitian connection on $L$
with curvature form
 \[
F_a=2\pi i\cdot\om,
 \]
 of Hodge type $(1,1)$ for the given holomorphic structure on $L$. Thus
the pair $(X,\om)$ is a symplectic manifold, the phase space of a
mechanical system. There are no invariants of an embedding of a Lagrangian
submanifold $\sL$ in a symplectic manifold. There are two ways of getting
invariants:
 \begin{enumerate}
 \item considering families of Lagrangian manifolds admitting invariants
(in particular limit singular subcycles); or
 \item giving submanifolds an additional structure (such as a section of
some bundle or an Hermitian connection on the trivial line bundle).
 \end{enumerate}
 The restriction of the pair $(L^k,a_{L^k})$ to a Lagrangian cycle $\sL$
gives this type of additional structure. It defines the space of covariant
constant sections:
 \[
 H^0_a ((L^k,a_{L^k})\rest{\sL})
 \]
which can be nonzero, as for points. Indeed, restricting to any Lagrangian
submanifold $\sL$ gives a topologically trivial line bundle on $\sL$ with
flat connection. A connection of this type is defined by its monodromy
character
 \begin{equation}
 \chi\colon \pi_1(\sL)\to \U(1),
 \label{eq:chi}
 \end{equation}
and admits a covariant constant section (as for restriction to a point) if
and only if this character is trivial.

 \begin{dfn} A Lagrangian cycle $\sL$ is a {\em level $k$
Bohr--Sommerfeld} (BS$_k$) cycle if the character (\ref{eq:chi}) for
$(L^k,a_{L^k})$ is trivial.
 \end{dfn}

In particular, just as for points,
 \[
\text{$\sL$ is BS$_k$} \implies H^0_a((L,a_L) \rest{\sL})=\C.
 \]

Moreover such a section defines a {\em trivialization} of the restriction 
$L^{k+d}\rest{\sL}$, which identifies $C^\infty$ sections with complex
valued functions on $\sL$:
 \[
 \Ga (L^{k+d}\rest{\sL})=C^{\infty}_\C(\sL).
 \]
Thus the restriction
to $\sL$ defines an embedding
 \begin{equation}
 \res\colon H^0(L^{k+d}) \hookrightarrow C^{\infty}_\C(\sL).
 \label{eq1.12}
 \end{equation}

 \begin{dfn} The image
 \begin{equation}
 \res(H^0(L^k))=\sH_\sL\subset C^\infty(\sL)
 \label{eq1.14}
 \end{equation}
is called the {\em analog of the Hardy space}.
 \end{dfn}

Now recall that our space $\sH_I^k$ (1.1) is the space of twisted
holomorphic {\em half-forms}:
 \[
 \sH^k_I=H^0(L^k\tensor L_{K/2}).
 \]
 To preserve the geometric meaning, fix a half-form $\hF$ on $\sL$; we
call a pair $(\sL,\hF)$ a {\em half-weighted Lagrangian} cycle or a
Lagrangian cycle marked with a half-form. Now we can identify the space of
functions with the space of half-forms
 \begin{equation}
 \Ga(L^{k+d})\rest{\sL}=C^{\infty}_\C(\sL)\cdot\hF=\Ga(\De^{1/2}),
 \label{eq:Ga}
 \end{equation}
where $\De$ is the complex volumes bundle on $\sL$. This space is
selfadjoint with respect to a Hermitian form like (1.8).

Following Borthwick, Paul and Uribe \cite{BPU}, we can construct a
distribution in some completion of $C^\infty(\sL)\cdot\hF$. Its
restriction to the image of $\sH_I^k$ gives a covector or a state. The BPU
method uses usual codimension 1 boundary conditions rather than Lagrangian
boundary conditions, and the original Hardy spaces of strictly
pseudoconvex domains rather than the analog of Hardy space (1.13). We
refer the reader to the beautiful paper \cite{BPU} for the details, which
we cannot reproduce here; this paper realizes a very large program. The
construction is following:
 \begin{enumerate}
 \item Our Hermitian connection on $L^*$ defines a contact structure on the
unit circle bundle $P$ of $L^*$.
 \item The disc bundle in $L^*$ is a strictly pseudoconvex domain, and
there is the Szeg\"o orthogonal projector $\Pi\colon L^2(P)\to \sH$ to the
Hardy space of boundary values of holomorphic functions on the disc bundle.
 \item The contact manifold $P$ is a principal $\U(1)$-bundle, and the
natural $\U(1)$-action on $P$ commutes with $\Pi$ and gives a
decomposition $\sH=\bigoplus_k H^0(L^{k+d})$ of the Hardy space.
 \item If we fix a metaplectic structure on $P$, we can lift every
$\BS^{k+d}$ submanifold to a Legendrian submanifold $\La\subset P$ over it,
marked with the lifted half-form $\hF$.
 \item $\La$ has an associated space of Legendrian distributions of order
$m$, which is the Szeg\"o projection of space of conormal distributions to
$\La$ of order $m + \half\dim X$ (see \cite{BPU}, 2.1).
 \item A half-form on $\La$ is identified with the symbol of a Legendrian
distribution of order $m$ (see \cite{BPU}, 2.2); thus at the level of
symbols, all Legendrian distributions look like delta functions or their
derivatives.
 \item For a Legendrian submanifold $\La$ with a half-form we fix the
Legendrian distribution of order $\frac{1}{2}$ with symbol $\hF$ which is
the Szeg\"o projection of the delta function $\de_\La$. 
 \end{enumerate}
In summary, we have:
 \begin{enumerate}
 \item For every lift $\La\subset P$ of a $\BS^{k+d}$ submanifold $\sL$
marked with a half-form $\hF$ we have a vector
 \begin{equation}
 \BPU_{k+d}(\La,\hF)=\Pi^k_{\hF} (\de_\La)\in H^0(L^{k+d}),
 \label{eq1.15}
 \end{equation}
where $\Pi^k_{\hF}$ is the Szeg\"o projection to the $(k+d)$th component of
the Hardy space of the distribution with symbol $\hF$.
 \item Every such lifting is defined up to $\U(1)$-action on $P$; thus a
pair $(\sL,\hF)$ defines a point of the projectivization
 \begin{equation}
 \BPU_{k+d} (\sL,\hF)=\PP (\Pi^k_{\hF} (\de_\La))\in \PP H^0(L^{k+d}).
 \label{eq1.16}
 \end{equation}
 \end{enumerate}
Our observations are the following:
 \begin{enumerate}
 \item This construction holds literally in the case that $\sL$ has the
structure of a smooth orbifold.
 \item If (1.3) holds, there exists a canonical {\em geodesic} lifting
(see Section~2). Thus $(\sL,\hF)$ defines a section of $L^{k+d}$. 
 \end{enumerate}

 The next step is the Analog of Serre's Theorems~A and~B proved in
\cite{BPU}, Section~3:

 \begin{thm} If\/ $k$ is large enough then\/ $\BPU_k(\La,\hF)\ne0$.
 \end{thm}

There are two or three canonical ways to give any Lagrangian submanifold
$\sL$ a half-form:
 \begin{enumerate}
 \item If $X$ is a K\"ahler manifold with a metaplectic structure. Then
this metaplectic structure defines a metalinear structure on $\sL$ (see
for example Guillemin \cite{Gu}), and the K\"ahler metric $g$ defines a
half-form $\hF_g$ on $\sL$. (This method is of course the most important
for our applications.)

 \item The graph of a metasymplectomorphism with symplectic volume as
square of the half-form.

 \item If $\sL$ admits a free torus action.
 \end{enumerate}

We have seen that half-weighted Bohr--Sommerfeld orbifolds look
geometrically like points. Let us denote by $\sL\sM$ the family of all
cycles that are Lagrangian with respect to $\om$. A polarization $(L,a_L)$
defines a subspace
 \begin{equation}
 \BS^{k+d}_l(L)\subset \sL\sM
 \label{eq:BS^{k+d}}
 \end{equation}
of Bohr--Sommerfeld Lagrangians; we decorate it by the index $l=[\sL]\in
H^{\dim_\C X}(X,\Z)$ (the cohomology class of the cycles), and by the level
$k$ (which we sometimes omit). This space breaks up into connected
components according to the topological type of generic Lagrangian cycles.
Recall that by the Darboux--Weinstein theorem we can identify a small
tubular neighborhood of $\sL$ with a neighborhood of the zero section of
the cotangent bundle of $\sL$, and any Lagrangian cycle in this
neighborhood can be identified with a {\em closed} 1-form on $\sL$. We get
a system of ``charts'' for $\sL\sM$, and the tangent space
 \begin{equation}
 T\sL\sM_\sL=\{\al\in \Om_\sL\bigm|\dd \al=0\}
 \label{eq1.17}
 \end{equation}
is the space of closed 1-forms on $\sL$. On the other hand, the periods of
these forms give infinitesimal deformations of the character (1.11).
Thus if $\sL\in\BS(L)$, the tangent space is the space of exact forms on
$\sL$:
 \begin{equation}
T\BS(X,L)_\sL=\{\al\in \Om_\sL\bigm|\al=\partial f\}=C^\infty (\sL)/\R.
\label{eq1.18}
 \end{equation}

 Thus we get the following result.

 \begin{prop} The normal space of\/ $\BS(L)$ in $\sL\sM$ at\/ $\sL$ is 
 \begin{equation}
N\BS(L)_\sL=H^1(\sL, \R).
\label{eq1.19}
 \end{equation}
 \end{prop}
 \begin{rmk} The subspace $\BS(L)$ of the space of all Lagrangian cycles
is a partial case of {\em isodrastic} deformations of Lagrangian cycles
(see \cite{Wei2}). \end{rmk}

Now for every family of Lagrangian cycles with base $B$, for a smooth
element $\sL$ we have the ``Kodaira--Spencer'' map
 \[
 \KS\colon TB_\sL\to H^1(\sL, \R).
 \]
The base $B$ of any family of Lagrangian cycles contains the subspace of BS
cycles
 \[
B\cap \BS^k(L)\subset B.
 \]
 \begin{cor} The codimension of this subset at a smooth cycle $\sL$ is 
 \[
\codim (B\cap \BS(L))=b_1 (\sL)-\corank \KS.
 \]
 \end{cor} 

We write
 \begin{equation}
\hWBS^{k+d}_l(L)
\label{eq1.20}
 \end{equation}
for the family of half-weighted $\BS$ cycles marked with half-forms, that
is, the set of pairs $\{(\sL,\hF)\}$. Then the BPU construction gives a
``rational'' map
 \begin{equation}
 \PP \fie_k\colon \hWBS_d^{k+d}(L)\to \PP H^0(L^{k+d})^*
 \label{eq1.21}
 \end{equation} 
which is regular if $k\gg0$ (just as the map (1.6) for points).

There are three types of finite dimensional families of Lagrangian cycles
where we can expect the existence of a finite set of BS cycles.

\paragraph{Example 1. Real polarization} A real polarization of
$(S,\om,L,a)$ is a fibration
 \begin{equation} \pi\colon S\to B,
\label{eq1.22}
 \end{equation}
such that $\om\rest{\pi\1(b)}=0$ for every point $b\in B$
and for generic $b$ the fibre $\pi\1(b)$ is a smooth Lagrangian.

Thus if we consider the pair $(S,\om)$ as the phase space of a mechanical
system, it admits a real polarization if and only if it is completely
integrable. In the compact case a generic fibre is a $n$-torus $T^n$
(where $2n=\dim_\R X$), and $\dim B=n$; thus
 \begin{equation}
\dim B=\rk H^1(T^n, \R) \implies \dim B\cap \BS(L)=0.
\label{eq1.23}
 \end{equation}
 \begin{rmk} A priori, there is no consistent way to introduce a preferred
orientation in the space of fibres of a real polarization. A metaplectic
structure provides it in some cases.
 \end{rmk}

\paragraph{Example 2. Moduli spaces of spLag cycles} Let $\sL$ be a
special Lagrangian cycle and $\sM^{[\sL]}$ the ``moduli space'' of all
deformation of $\sL$ as a special Lagrangian cycle in $X$ (see \cite{T3}).
Then by McLean's theorem the tangent space $T \sM^{[\sL]}=H^1(\sL,\R)$ is
the space of harmonic 1-forms and the Kodaira--Spencer map has
$\corank\KS=0$. Therefore, by definition, every smooth BS cycle must be 
{\em infinitesimally rigid}, so that
 \begin{equation}
\dim (\sM^{[\sL]}\cap \BS_l)(L)=0.
 \label{eq1.24}
 \end{equation}

In particular, if $X$ is a Calabi--Yau threefold polarized by a Ricci flat
metric with fixed complex orientation we have the system of functions of any level $k$
 \begin{equation}
H^3(X, \Z)\to \Z
\label{eq1.25}
 \end{equation}
sending a cohomology class $l\in H^3(X, \Z)$ to the number 
 \[
 \#(\sM^d\cap \BS_l^k(L))
 \]
(see \cite{T4} and \cite{T2} for the relation of this function with the
Casson--Donaldson invariant).

Now fixing a half-form on all the $\BS^{k+d}$ cycles of these families, we
get a finite set of points in $\PP H^0(L^{k+d})$.

Moreover in many cases this collection of points in $\PP H^0(L^{k+d})$ can
be lifted up to finite ambiguity to a basis of the vector space
$H^0(L^{k+d})$ (as predicted in \cite{BPU}, Remark on p.~400).

It was realized by Poincar\'e for the case $X=C$ is an algebraic curve of
genus $g > 1$.

\paragraph{Example 3. Relative Poincar\'e series} Let $C$ be an
algebraic curve of genus $g>1$ with a fixed $\Spin^\C$ structure, that is,
with a fixed {\em theta structure} $L$ such that
 \[
 L^2=L_{K_C}.
 \]
Then $L$ defines a metaplectic structure and a polarization.

This line bundle has a Hermitian connection with square the Levi-Civita
Hermitian connection on $T^*C$. Then every 1-cycle $\sL$ is Lagrangian.
 \[
\sL\in \BS^2(L) \quad \text{if and only if it is
{\em geodesic}}.
 \]
 If $k>1$ then it is $\BS^k$ if and only if it is {\em $k$-geodesic}, that
yis, its holonomy is a $k$th root of unit.

 Parametrizing such a cycle by arclength gives a half-form $\hF$ on it.
Thus it defines the Poincar\'e series of $H^0(L^{2k})$ as an automorphic
form given by the relative Poincar\'e series (see \cite{BPU}, Section~4).

In the same vein, a $\BS^{k+d}$ cycle $\sL$ defines a section of $L^{k+d}$
if the canonical class $[K_X]=d c_1(L)$ with $d\in \Z$ and the Hermitian
connection $a$ is ``proportional'' to the connection on $\det T^*X$ induced
by the Levi-Civita connection (see below).

 \begin{rmk} As in the original proof of Serre's Theorems~A and~B, we
can't avoid some technical work in functional analysis. Our aim is to
localize these techniques in one place, the proof of Theorem~1.1. After
this the theory becomes a combination of projective algebraic geometry and
Lagrangian geometry. We call this hybrid {\em aLag geometry}.
 \end{rmk}

\section{BPU construction, geodesic lifting and\\ geo\-metric quantization}

For our applications, we extend slightly the BPU construction described in
(1.15--16) of the previous section. We must repeat some of the details. We
stay in the situation of the starting point of Section~1: let $X$ be a
K\"ahler manifold with a polarization $L$ having a Hermitian connection
with the K\"ahler form as curvature form.
 
Consider the principal $\U(1)$-bundle $P$ of the dual line bundle $L^*$,
the unit circle bundle in $L^*$. Let
 \[
D\subset L^*, \quad \partial D=P
 \]
 be the unit disc subfibration with boundary. Our $L$ is positive,
hence $D$ is a strictly pseudoconvex domain. The Hermitian connection on
$L^*$ is given by 1-form $\al$ on $P$ which defines a contact structure on
$P$ with volume form
 \[
 \frac{1}{2\pi} \al \wedge \dd\al^n, \quad \text{where} \quad n=\dim_{\C}
X.
 \]
The null space of $\al$ at a point $p\in P$ is the maximal complex
subspace of the tangent space.

 The Hardy subspace
 \begin{equation}
 \sH_I\subset L^2(P)
 \label{eq2.1}
 \end{equation}
consists of boundary values of holomorphic functions on $D$. We have
the Szeg\"o orthogonal projector
 \begin{equation}
 \Pi_I\colon L^2(P)\to \sH_I.
 \label{eq2.2}
 \end{equation}
The natural action of $\U(1)=\Aut P$ on $P$ as a principal bundle commutes
with $\Pi_I$ and decomposes the space $\sH_I$ as a Hilbert direct sum of
isotypes:
 \begin{equation}
\sH_I=\bigoplus_{k=0}^{\infty} (\sH_I^k=H^0(L^{k+d})),
 \label{2.3}
 \end{equation}
with only positive characters. Thus
 \begin{equation}
 \Pi_I=\bigoplus_{k=0}^{\infty} \Pi_k.
 \label{eq2.4}
 \end{equation}
Our first addition to \cite{BPU} is the following: suppose that the group
$G$ acts on $X$ preserving $\om$ and $(L,a_L)$. Then there exists a central
extension $\wG$:
 \begin{equation}
 1\to \U(1)\to\wG\to G\to 1,
 \label{eq2.5}
 \end{equation}
where the centre $\U(1)=\Aut P$ acts on $P$ as a group of a principal
bundle. This action induces a natural representation
 \begin{equation}
 \rho\colon\wG\to\Op(L^2(P))
 \label{eq2.6}
 \end{equation}
on the operator algebra of the space of functions.

If the transformations of $G$ preserve our complex structure $I$ then the
projector $\Pi_I$ (2.4) defines a representation
 \begin{equation}
 \Pi_I \circ \rho \circ \Pi_I\colon\wG\to\Op(\sH_I)
 \label{eq2.7}
 \end{equation}
which commutes with the action of $\U(1)$. So this action decomposes this
representation as a Hilbert direct sum of isotypes
 \begin{equation}
 \rho_I^k=\Pi_k \circ \rho \circ \Pi_k\colon\wG\to
 \End H^0(L^k).
 \label{eq2.8}
 \end{equation}
These representations can be projectivized 
 \[
 \PP\rho_I^k\colon G\to\Aut\PP H^0(L^k).
 \] 
Recall that our $X$ and the principal bundle $P$ have given metaplectic
structures.

Now if $\sL$ is a half-weighted BS orbifold, the complex conjugate of a
covariant constant section gives a lift of it to a Legendrian cycle $\La$
on $P$ marked with half-forms $(\La,\hF)$, and the construction (1.15)
defines a section
 \[
 \BPU_{1-d}(\La,\hF)\in H^0(X, L).
 \]
 In the same vein, we get lifts $\La_{k+d}$ of $\sL$ from $\BS^{k+d}(L)$
to the Legendrian cycle on $P$ which is a cyclic $(k+d)$-cover of it and
the system of sections
 \begin{equation}
\BPU_{k+d}(\La_{k+d},\hF)\in H^0(X, L^{k+d}).
 \label{eq2.9}
 \end{equation}

 The Lagrangian cycle $\sL$ can be reconstructed from the set of its BPU
images as the quasi-classic limit as $1/k=\text{Planck's constant}\to0$:
the wave fronts of distributions concentrate on $\sL$ (see \cite{BPU} and
the references given there).

Now for two Lagrangian cycles $(\sL_1,\hF_1)$ and $(\sL_2,\hF_2)$, the
asymptotic behaviour of the scalar product $\Span{\BPU_k(\sL_1,\hF_1),
\BPU_k(\sL_2,\hF_2)}$ (see (1.8)) can be computed in terms of the
intersection $\sL_1\cap \sL_2$ (see \cite{BPU}). In particular,
 \begin{equation}
\sL_1\cap \sL_2=\emptyset \implies \BPU_k(\sL_1,\hF_1) \perp
\BPU_k(\sL_2,\hF_2)
\label{eq2.10}
 \end{equation}
asymptotically as $k\to\infty$. (For the orbifold case these asymptotics
are somewhat weaker, but are still quite expressive for geometric
corollaries). This asymptotic technique comes from the physical
interpretation of this set-up as the ``classical'' (= pre-BRST) geometric
quantization (GQ for short). More precisely the asymptotic analysis of
quantum states gives
 \begin{multline}
 \Span{\BPU_k(\sL,\hF), \BPU_k(\sL,\hF)}^k \ \sim\ 
 \left(\frac{k}{\pi}\right)^{\half\dim X}\int_{\sL} \Vert\hF \Vert^2 \\
 +O(k^{\half(\dim X-1)}),
 \label{eq2.11}
 \end{multline}
 and if $\sL_1\cap \sL_2=\emptyset$ then
 \begin{multline}
 \Span{\BPU_k(\sL_1,\hF_1), \BPU_k(\sL_2,\hF_2)}^k \ \sim\ 
 \left(\frac{k}{\pi}\right)^{\half(\dim X-1)}+ \\
 O(k^{\half(\dim X-2)}).
 \label{eq2.12}
 \end{multline}
 From this it easy to see that
 \begin{prop} Let $\pi\colon X\to B$ be any real polarization (see
Section~1, Example~1). Then for $k\gg0$, the BPU vectors in $H^0(L^{k+d})$
span a subspace of dimension
 \[
 \rk \Span{\BPU_{k+d}(\BS^{k+d}(L)\cap B)}=\#(\BS^k(L)\cap B).
 \]
 \end{prop}

 \begin{rmk} A more sophisticated analysis of the asymptotics of quantum
states extends this observation as follows: let
 \[
\sL_1,\dots, \sL_{N_{\max}}\subset\BS^{k+d}(L)
 \]
be a maximal collection of {\em disjoint\/} $\BS^{k+d}$ cycles. Then
 \[
N_{\max} \le\rk H^0(L^{k+d}),
 \]
and the right-hand side is given by the Riemann--Roch theorem. \end{rmk}

\subsection*{Geodesic lifting}
The lifting of a BS Lagrangian cycle on $X$ to a Legendrian cycle on $R$
we have described is defined up to the natural $\U(1)$-action on $R$ and
the states $\BPU_k(\sL,\hF)$ are defined up to scaling. But in our
applications we can do this almost canonically (up to a finite ambiguity)
and get an actual basis of $H^0(L^{k+d})$.

To describe this almost canonical lifting we must consider the Lagrangian
Grassmannization of the tangent bundle of $X$ as described in \cite{T3}.
Pointwise, the tangent space $(TX)_x$ at a point $x\to X$ is
$\C^n$ with the constant symplectic form $\Span{\ \,,\ }=\om_x$ and the
constant Euclidean metric $g_x$, giving the Hermitian triple
$(\om_x,I_x,g_x)$. Define the Lagrangian Grassmannian $(\LaGr)_x
=\LaGr(TX)_x$ to be the Grassmannian of oriented Lagrangian
subspaces in $(TX)_x$. Taking this space over every point of
$X$ gives the oriented Lagrangian Grassmannization of
$TX$
 \begin{equation}
 \pi\colon \LaGr(TX)\to X
\quad\text{with} \quad \pi\1 (x) =(\LaGr)_x.
 \label{eq2.13}
 \end{equation}
A complex structure $I_x$ on $(TX)_x$ gives the standard
identification
 \begin{equation}
(\LaGr)_x=\U(n) / \SO(n).
 \label{eq2.14}
 \end{equation}
This space admits a canonical map
 \begin{equation}
 \det\colon(\LaGr)_x\to \U(1)=S^1_x
 \quad\text{sending $u\in\U(n)$ to $\det u\in\U(1)=S^1$.}
 \label{eq2.15}
 \end{equation}
 Taking this map over every point
of $X$ gives the map
 \begin{equation}
 \det\colon \LaGr(TX)\to S^1(L_{-K}),
 \label{eq2.16}
 \end{equation}
where $S^1(L_{-K})$ is the unit circle bundle of the line bundle
$\bigwedge^{n} TX =\det TX$, with first Chern class
 \[
 c_1(\det TX)=-K_X,
 \]
where $K_X$ is the canonical class of $X$ (see for example
\cite{T2} and \cite{T3}).

We have already noted that our Lagrangian cycles does not usually have an
orientation defined a priori. Thus we must consider the Lagrangian
Grassmannian $\La(TX)$ forgetting orientations. Then we get a map
 \[
\det\colon \La(TX)\to S^1(L_{-K/2}) 
 \]
in place of (2.16).

Now for every oriented Lagrangian cycle $\sL\subset X $, we have the
Gaussian lift of the embedding $i\colon \sL\to X$ to a section
 \begin{equation}
 G(i)\colon \sL\to \La(TX)\rest{\sL},
 \label{eq2.17}
 \end{equation}
 sending $x\in \sL$ to the subspace $T\sL_x\subset(TX)_x$. The composite
of this Gauss map with the projection $\det$ gives the map
 \begin{equation}
 {\det}\circ G(i)\colon \sL\to S^1(L_{-K/2})\rest{\sL}.
 \label{eq2.18}
 \end{equation}
Thus every Lagrangian cycle $\sL$ defines a Legendrian subcycle
 \begin{equation}
 \La={\det}\circ G(i) (\sL)\subset S^1(L_{-K/2}).
\label{eq2.19}
 \end{equation}
The Levi-Civita connection of the K\"ahler metric defines a Hermitian
connection $a_{\LC}$ on $L_{-K/2}$.

 \begin{dfn} A Lagrangian cycle $\sL$ is {\em almost geodesic} if the
Legendrian cycle $\La={\det}\circ G(i)(\sL)$ is horizontal with
respect to the Levi-Civita connection $a_{\LC}$ on $L_{-K/2}$. \end{dfn}

We now use property (1.3). The line bundle $L_{-K/2}$ is $L^{-d}$, where
$L$ is the line bundle of the polarization and we suppose that the
Levi-Civita connection is induced by the connection $a_L$ on $L$. Then we
have
 \begin{prop} A Lagrangian cycle $\sL$ is $\BS^0$ if and only if it is
almost geodesic.
 \end{prop}

 The Hermitian structures of our line bundles define a map
 \begin{equation}
 \mu_d\colon S^1(L^*)\to S^1(L_{-K/2})
 \label{eq2.20}
 \end{equation}
of the principal $\U(1)$-bundles of these line bundles, which fibrewise is
minus the isogeny of degree $d$. Thus {\em every} Lagrangian cycle $\sL$
defines an oriented Legendrian subcycle (see (2.1))
 \begin{equation}
 \La=\mu_d\1 ({\det}\circ
G(i) (\sL))\subset S^1(L^*)=P.
 \label{eq2.21} \end{equation}
 
Now consider the pair of isogenies
 \begin{equation}
\mu_{d+k}\colon S^1(L^d)\to S^1(L^{d(k+d)})
 \label{eq2.22}
 \end{equation}
and
 \[
\mu_d\colon S^1(L^{k+d})\to S^1(L^{d(k+d)})
 \]
and the lift 
 \[
l\colon \sL\to S^1(L^{k+d})
 \]
given by a covariant constant section over a $\BS^{k+d}$ cycle $\sL$.
 \begin{dfn} The lift $l$ is {\em almost geodesic} if
 \[
 \mu_{d+k} \circ {\det}\circ G(i) (\sL)=\mu_d \circ l (\sL).
 \]
 \end{dfn}

The number of geodesic lifts is obviously $\le|d(k+d)|$.

In summary, let $\La^k\sM$ be the space of Legendrian subcycles of $P$
the images of whose projection to $X$ is the $k$th root of unity cover of
a $\BS^{k+d}$ Lagrangian cycle on $X$ (such Legendrian cycles are
sometimes called {\em Planckian cycles}). Then the natural projection
 \[
p\colon \La^{k+d}\sM\to \BS^{k+d}(L)
 \]
which sends a Legendrian cycle to Lagrangian cycle is a principal
$\U(1)$-bundle. 

The geodesic lifting 
 \begin{equation}
l\colon \wBS^{k+d}(L)\to \La^{k+d}\sM
\label{eq2.23}
 \end{equation}
we have described is a multisection of this principal bundle and 
 \[
p\colon \wBS^{k+d}(L)\to \BS^{k+d}(L) 
 \]
is a finite cyclic cover.

Consider a real polarization $\pi\colon X\to B$ (1.23). Then we have a finite
set of Bohr--Sommerfeld fibres
 \[
 B\cap \BS^{k+d}(L)=\{\sL_i\}, \quad \text{for $i=1,\dots, N_\pi^{k+d}$.}
 \]

 \begin{dfn} A choice of geodesic lifts 
 \[
\{\wsL_i\}\subset \La^{k+d}\sM 
 \]
 is called a choice of {\em theta structure} of the real polarization
$\pi$.
 \end{dfn}

Marking these Lagrangian cycles with the half-forms given by our K\"ahler
metric $g$ (see (1) below Theorem 1.1), we get a finite set 
 \[
 \{\wsL_i,\hF_g\}
 \]
of half-weighted Legendrian cycles.

We know that 
 \[
\BPU_k (\{\wsL_i,\hF_g\})\subset \PP H^0(L^{k+d})
 \]
is a linear independent system of vectors (states) if $k\gg0$. In
particular, if
 \begin{equation}
 \#(\BS^{k+d}(L)\cap B)=\rk H^0(L^{k+d}),
 \label{eq2.24} \end{equation}
 we get a Bohr--Sommerfeld basis.

 \begin{rmk} It is easy to see that we are imitating the geometric
situation of Section~1, Example~3. For other descriptions and applications
of the geodesic lift from Lagrangian to Legendrian cycles see \cite{T1},
\cite{T2} and \cite{T3}.
 \end{rmk}
 
\subsection*{Geometric quantization}

There is a deep reason for coincidences such as (2.24) for
Bohr--Sommerfeld fibres of a real polarization of the phase space of a
classical mechanical system: we can view any symplectic manifold $(S,\om)$
as the phase space of some classical mechanical system, and the pair
$(L,a_L)$, where $a_L$ is an Hermitian connection on line bundle $L$ with
curvature form $F_a=2\pi i\cdot\om$ as a prequantization data of this
system.

Bohr--Sommerfeld bases identify two approaches to the geometric
quantization of $(S,\om,L,a)$ (see \cite{A}, \cite{S1} or \cite{W}). The
first approach is a choice of a complex polarization, which is nothing
other than a choice of a complex structure $I$ on $S$ such that $S_I=X$ is
a K\"ahler manifold with K\"ahler form $\om$. Then the curvature form of
the Hermitian connection $a$ is of type $(1,1)$, hence for any level
$k\in\Z^+$, the line bundle $L^k$ is a holomorphic line bundle on $S_I$.
Complex quantization provides the space of wave functions of level $k$
(1.1) (see Kirillov's survey \cite{K}).

The second approach to geometric quantization is the choice of a real
polarization of $(S,\om,L,a)$ (see Section~1, Example~1) which is a
fibration $\pi\colon S\to B$ (1.23). We have already seen that restricting
$(L,a)$ to a Lagrangian fibre gives a flat connection or equivalently, a
character of the fundamental group $\chi\colon\pi_1(\text{fibre})\to\U(1)$.

Let $\sL_{\pi}$ be the sheaf of sections of $L$ that are covariant
constant along fibres. Then we get the space $\sH_{\pi}=\bigoplus_{i}
H^i(S,\sL_{\pi})$ and in the regular case, \'Sniatycki proved that
$H^i(S,\sL_{\pi})=0$ for $i\ne n$. To compute the last component
$H^n(S,\sL_{\pi})$ we need to involve Bohr--Sommerfeld fibres: we expect to
get a finite number of Bohr--Sommerfeld fibres, and in the regular case,
 \[
 H^n(S,\sL_{\pi})=\bigoplus_{{\BS}\cap B}\C\cdot s_{i},
 \]
where $s_{i}$ is the covariant constant section of the restriction of
$(L, a)$ to a Bohr--Sommerfeld fibre of the real polarization $\pi$ (see
\cite{S2}).

In the general case, we can use this to define a new collection of
spaces of wave functions (of level $k$):
 \begin{equation}
 \sH_{\pi}^k=\bigoplus_{{\BS^{k+d}}\cap B}\C\cdot s_{i},
 \label{eq2.25}
 \end{equation}
and use the Borthwick--Paul--Uribe construction to compare (1.1) with
(\ref{eq2.25}).

There is a canonical way of describing the subset $\BS^{k+d}(L)\cap B$
using special coordinates on $B$, the so-called {\em action} coordinates,
which are part of the action angle coordinates (see \cite{A},
\cite{GS1},\dots).

 An important observation, proved mathematically in some cases, is that
the projectivization of the spaces (\ref{eq1.1}) does not depend on the
choice of complex structure:
 \begin{equation}
 \frac{\partial\PP\sH_{I}^k}{\partial I}=0.
 \label{eq2.26}
 \end{equation}
In other words, the vector bundle (1.2) admits a projective flat
connection. Thus spaces of wave functions are given purely by the
symplectic prequantization data (see for example \cite{H}). The same is
true for the projectivization of the spaces (\ref{eq2.25}). Moreover,
these spaces do not depend on the real polarization $\pi$ (1.21), provided
that we extend our prequantization data $(S,\om,L,a,)$ by adding some
half-density or half-form on every BS fibre (see \cite{GS1}) to define the
half-form pairing of Blattner, Kostant and Sternberg (for the difference
between half-density and half-form quantizations see \cite{W}).

To compare the spaces
 \[
 \sH_I^k \quad \text{and} \quad \sH_{\pi}^k
 \]
by the BPU method, we have to arrange for them to have the same rank. This
arithmetical problem can be solved directly in many interesting cases;
see, for example, \cite{K1} and \cite{JW1}. For the geometry behind these
coincidences see \cite{T1}.

 \section{Application: theory of non-Abelian theta functions}

 The classical theory of theta functions serves as a beautiful model for
our theory. Although this theory is realized by many approaches to
geometric quantization (see \cite{T1}), we must demonstrate that all
classical bases of theta functions can be described as Bohr--Sommerfeld
bases given by the BPU method, that the geodesic lifting is a choice of
theta structure and so on. This is a beautiful but quite serious job, and
will be done in a special paper (or book). Here we will discuss this
theory as a model for the theory of {\em non-Abelian theta functions}.

 Let $A$ be a principal polarized Abelian variety (ppAv) of complex
dimension $g$ with zero element $o\in A$ and with flat metric $g$. Then
the tangent bundle $TA$ has the standard constant Hermitian structure
(that is, the Euclidean metric, symplectic form and complex structure
$I$). The K\"ahler form $\om$ gives a polarization of degree 1. In the
equality (1.3) we have $d=0$. We fix a smooth Lagrangian decomposition of
$A$
 \begin{equation}
 A=T^g_+\times T^g_-, \label{eq3.1}
 \end{equation}
 with both tori Lagrangian with respect to $\om$ (recall that, smoothly,
$A$ is the standard torus $\R^{2g} / \Z^{2g}$ with the standard constant
integer form $\om$ and this decomposition is nothing other than reducing
the integer form $\om$ to normal form). Let $L$ be a holomorphic line
bundle with holomorphic structure given by a Hermitian connection $a$ with
curvature form $F_a=2\pi i\cdot\om$, and $L=\Oh_A(\Theta)$, where $\Theta$
is the classical symmetric theta divisor. The decomposition
(\ref{eq3.1}) induces a decomposition $H^1(A,\Z)=\Z^g_+\times\Z^g_-,$ and
a Lagrangian decomposition
 \begin{equation} A_k=(T^g_+)_k \times (T^g_-)_k
 \label{eq3.2} \end{equation}
 of the $k$-torsion subgroup.

In this case, complex quantization is nothing other than the classical
theory of theta functions. Indeed, the Lagrangian decomposition
(\ref{eq3.2}) of the $k$-torsion subgroup defines a collection of
compatible theta structures of every level $k$ and a decomposition of the
spaces of wave functions
 \begin{equation}
 \sH_I^k=H^0(A, L^k)=\bigoplus_{w\in (\Z^g)_k^-}\C\cdot\theta_w
 \quad\text{with}\quad \rk\sH_I^k=k^g,
 \label{eq3.3}
 \end{equation}
where $\theta_w$ is the theta function with characteristic $w$ (see
\cite{Mum}).

On the other hand, the direct product (\ref{eq3.1}) gives us a real
polarization
 \begin{equation}
\pi\colon A\to T^g_-=B.
 \label{eq3.4}
 \end{equation}
Remark that in this case the action coordinates are just flat coordinates
on $T^g_-=B$, and under this identification
 \begin{equation}
 B\cap \BS^k(L)=(T^g_-)_k
 \label{eq3.5}
 \end{equation}
is the $k$-torsion subgroup.

Thus applying geometric quantization to the real polarization
(\ref{eq3.4}) of the phase space $(A,\om, L^k,a_{L^k})$, where $a_k$ is the
Hermitian connection defining the holomorphic structure on $L^k$, we get
the decomposition
 \begin{equation}
 \sH_{\pi}^k=\bigoplus_{\rho\in\U(1)^g_k}\C\cdot s_{\rho}.
 \label{eq3.6}
 \end{equation}

 \begin{cor} \begin{enumerate} \item $\rk \sH_{L^k}=\rk \sH_{\pi}^k=k^g$.

 \item Moreover, there exists a isomorphism
 \[
 \sH_I^k =\sH_{\pi}^k,
 \]
and this is canonical up to a scaling factor.
 \end{enumerate}
 \end{cor}
Indeed, the identification of the $\BS^k$ fibres of $\pi$ given by the
projection to $T^g_+$ gives us at the same time a lift of the $\BS^k$
fibres to Legendrian submanifolds of $P$. Moreover the canonical class
$K_A=0$. So there exists a canonical metaplectic structure on $A$, and a
canonical collection of half-forms on Bohr--Sommerfeld fibres invariant
with respect to translations defined up to a common phase factor. We can
use the Borthwick--Paul--Uribe homo\-morphism (2.23) which is an inclusion
 \begin{equation}
\BPU_k\colon \sH_{\pi}^k \hookrightarrow \sH_I^k
\label{eq3.7}
 \end{equation}
(because sections are orthogonal) and which is an isomorphism (because the
ranks are equal). Moreover it easy to prove

 \begin{prop} The homomorphism $\BPU_k$ (\ref{eq3.7}) extends to an
inclusion of\/ $H_k$-modules, where $H_k$ is the Heisenberg group of level
$k$.
 \end{prop}
 \begin{cor} For every $\BS_k$ torus $T^g_j$ of (\ref{eq3.5}) the BPU
quantum mechanical state $\BPU_k(T^g_j)$ coincides with the corresponding
theta function.
 \end{cor}

The functions making up the special bases of these spaces are called
classical theta functions of level $k$ with characteristic.

 \begin{rmk}
 This coincidence should of course be proved directly using the form of
the Schwartz kernel of $\Pi_k$, Fourier images of delta functions and the
interpretation of theta functions as solutions of the heat equation.
 \end{rmk}

Thus combining the constructions of complex and real polarizations gives
us some orthogonal bases in complete linear systems. However, if we start
with any polarized K\"ahler manifold $X$, the main question is about a real
polarization of the form (1.23) on $X$ (possibly with degenerate fibres).

Finally, let $\Si$ be a Riemann surface of genus $g$ and $A=J_{\Si}$ its
Jacobian. Then as a real manifold
 \[
J_{\Si}=T^{2g}=\Hom(\pi_1(\Si),\U(1)),
 \]
with the symplectic form $\om$ and the line bundle $L$ with a Hermitian
connection $a$ with curvature $F_a=2\pi i \om$. Thus we can apply
these constructions to the quadruple
 \begin{equation}
(\Hom (\pi_1(\Si),\U(1)), \om, L_\Theta, a).
\label{eq3.8}
 \end{equation}
 Then the Lagrangian decomposition (3.1) gives a real polarization with
$\BS^k$ fibres (3.5).

Now giving $\Si$ a complex structure $I$ defines a complex polarization of
$\Hom (\pi_1(\Si),\U(1))=J_{\Si}$. So the collection of spaces (1.1)
 \[
\sH^k_I=H^0(J_{\Si}, L^k_\Theta)
 \]
are fibres of the holomorphic vector bundles (1.2)
 \begin{equation}
\sH^k\to \sM_g
\label{eq3.9}
 \end{equation}
over the moduli space of Riemann surfaces of genus $g$.

Then the identification (3.7) shows that these spaces are actually
independent of the complex structure $I$; that is, there exists a
projective flat connection on every vector bundle (3.8). These connections
may be described by a heat equation as in \cite{H} and \cite{We}.

Here we apply this method to the following noncommutative generalization
of this situation: consider the $(6g-6)$-manifold
 \begin{equation}
 R_g=\Hom (\pi_1(\Si_2),\SU(2)) /\PU(2),
 \label{eq3.10}
 \end{equation}
the space of classes of $\SU(2)$-representations of the fundamental group
of this Riemann surface (it only depends on $g$) and apply the BPU
construction to its GQ.

The new feature of this situation is the fact that $R_g$ is singular. So
before this case, we must consider as a model the following {\em singular
Abelian case}. Let
 \begin{equation}
K_\Si=J_\Si / \{\pm \id\}
 \label{eq3.11}
 \end{equation}
be the Kummer variety. All geometric objects (3.8) are invariant with respect
 to the involution $-\id$ and we get the prequantized mechanical system
 \begin{equation}
(K_\Si, \om, L_\Theta, a)
 \label{eq3.12}
 \end{equation}
with singular phase space 
 \begin{equation}
 \Sing K_\Si=(J_\Si)_2,
 \label{eq3.13}
 \end{equation}
that is, the 2-torsion points of the Jacobian.

A complex structure $I$ on $\Si$ defines a complex polarization of
$K_\Si$; but now we only consider spaces of wave functions of even level
 \begin{equation}
 \sH_I^{2k}=H^0(K_\Si, L^{2k}_\Theta)=H^0_\ev(J_\Si,L^{2k}_\Theta),
 \label{eq3.14}
 \end{equation}
the space of even (symmetric) theta functions.

To describe a real polarization of $K_\Si$, represent $\Si$ as a
connected sum of $g$ 2-tori:
 \begin{equation}
 \Si=T^2_1 \cosum T^2_2 \cosum \cdots\cosum T^2_g,
 \label{eq3.15}
 \end{equation}
and fix a standard pair of generators of the fundamental group $(a_i,b_i)$
of each 2-torus $T^2_i$. Then we get the standard presentation of the
fundamental group of $\Si$
 \begin{equation}
 \pi_1(\Si)=\Span{a_1,\dots, a_g, b_1,\dots,
b_g\bigm|\prod\nolimits_{i=1}^{g} [a_i, b_i]=\id}.
 \label{eq3.16}
 \end{equation}
Killing the generators $a_i$ defines the {\em handlebody} $\wSi_a$ with
boundary
 \begin{equation}
 \partial \wSi_a=\Si,
 \label{eq3.17}
 \end{equation}
and fundamental group
 \begin{equation}
 \pi_1(\wSi_a)=\Span{b_1,\dots, b_g},
 \label{eq3.18}
 \end{equation}
 the free group on the $b_i$.

Now we can define the Jacobian of a handlebody
 \begin{equation}
 J_{\wSi_a}=H_1(\wSi_a, \R) / H_1(\wSi_a, \Z) 
 \label{eq3.19}
 \end{equation}
and its Kummer variety:
 \begin{equation}
 K_{\wSi_a}=J_{\wSi_a} / \{\pm\id\}.
 \label{eq3.20}
 \end{equation}
Our real polarization (3.4) can be described as the natural map
 \begin{equation}
 \pi\colon J_\Si\to J_{\wSi_a}
 \label{eq3.21}
 \end{equation}
providing a real polarization of the Kummer variety
 \begin{equation}
 \pi\colon K_\Si\to K_{\wSi_a}.
 \label{eq3.22}
 \end{equation}
The fibres of this polarization are Lagrangian, and the $2^g$ Kummer
varieties of $g$-dimensional tori over $(K_{\wSi_a})_2$ which are
singular, and
 \begin{equation}
 \Sing \pi\1 (w)=(\pi\1(w))_2, \quad\text{for $w\in (J_{\wSi_a})_2$.}
 \label{eq3.23}
 \end{equation}
Obviously the involution $-\id$ on $J_{\wSi_a}$ preserves the BS
fibres of (3.4) and acts freely on 
 \begin{equation}
 (J_{\wSi_a})_{2k}\setminus(J_{\wSi_a})_{2} 
 \label{eq3.24}
 \end{equation}
preserving pointwise the subset
 \[
 (J_{\wSi_a})_{2}.
 \]
 So the number of $\BS_{2k}$ fibres of the real polarization (3.21)
 \begin{equation}
 \#(K_{\wSi_a}\cap \BS^{2k}(K_\Si, L_\Theta))=2^{g-1}(k^g + 1)
 \label{eq3.25}
 \end{equation}
is equal to rank of the space of even theta functions of level $2k$ (3.14).

Now the BS fibres over 
 \begin{equation}
 \Bigl((J_{\wSi_a})_{2k}\setminus(J_{\wSi_a})_{2}\Bigr)\bigm/\{\pm\id\}
 \label{eq3.26}
 \end{equation}
are nonsingular, and the fibres over $(J_{\wSi_a})_{2}$ are singular
and {\em simply connected}. The following statement also holds in the
non-Abelian case:
 \begin{prop} For a singular $\BS_k$ fibre $\pi\1(w)$,
 \[
H_1(\pi\1(w))=H_1(\pi\1(w)\cap \Sing K_\Si). 
 \]
 \end{prop}
As usual we have the space of wave functions (like (3.6)) for the real
polarization (3.22)
 \begin{equation}
\sH^{2k}_\pi=\bigoplus_{w\in ((J_{\wSi_a})_{2k}\setminus(J_{\wSi_a})_{2}) /
 \{\pm\id\}}\C \cdot s_w \bigoplus_{w\in (J_{\wSi_a})_2} \C \cdot s_w.
 \label{eq3.27}
 \end{equation}

Summarizing, for even level we get the following orthogonal decomposition
of the space (3.3) of theta functions
 \begin{equation}
\sH_I^{2k}=H^0_\ev(J_\Si, L^{2k}) \oplus H^0_\odd(J_\Si, L^{2k}) 
 \label{eq3.28}
 \end{equation}
into even and odd theta functions, and the even component is the space of
wave functions (3.14) of a complex polarization of the Kummer variety.
Then the direct BPU construction gives a linear embedding
 \begin{equation}
BPU\colon \bigoplus_{w\in ((J_{\wSi_a})_{2k}\setminus(J_{\wSi_a})_{2}) /
 \{\pm\id\}}\C \cdot s_w\to H^0_\ev(J_\Si, L^{2k}_\Theta) 
 \label{eq3.29}
 \end{equation} 
for nonsingular BS fibres and a slight modification of it for the {\em
orbifold} case gives us the full identification
 \begin{equation}
H^0(K_\Si, L^{2k}_\Theta)=\sH^{2k}_\pi. 
\label{eq3.30}
 \end{equation} 
Under this identification, (3.27) corresponds to the decomposition
 \[
 H^0(L^{2k}_\Theta)=H^0(\Oh_{(J_\Si)_2})\oplus H^0 (J_{(J_\Si)_2}\tensor
L^k).
 \]

Unfortunately in the non-Abelian case the singularities are much worse than
in the case of orbifolds and we apply the following strategy: for the
smooth and orbifold cases we use the BPU method directly and its orbifold
modification, but for heavy singularities we use the special features of
our situation avoiding analysis.

 \subsection*{Complex quantization of $R_g$}

The space (3.10) is stratified by the
subspace of reducible representations
 \begin{equation}
 R^\triv\subset R_g^\red\subset R_g, \quad R_g\irr=R_g-R_g^\red.
 \label{eq3.31}
 \end{equation}
 Using symplectic reduction arguments, we get a nondegenerate closed
symplectic form $\Om$ on this space. This form $\Om$ defines a symplectic
structure on $R_g\irr$.

There exists a Hermitian line bundle $L$ with the $\U(1)$-connection
$A_{\CS}$ on $R_g$ (the Chern--Simons connection, see \cite{RSW} or
\cite{T1}, \S3). By definition, the curvature form of this connection is
 \begin{equation}
 F_{A_{\CS}}= 2 \pi i \cdot\Om.
 \label{eq3.32}
 \end{equation}
Thus the quadruple
 \begin{equation}
 (R_g,\Om, L, A_{\CS})
 \label{eq3.33}
 \end{equation}
is a prequantum system.

The standard way of getting a complex polarization is to give the Riemann
surface $\Si$ of genus $g$ a conformal structure $I$. We get a complex
structure on the space of classes of representations $R_g$ such that
$R_{\Si}=R_g=\sM\semis$ is the moduli space of semistable holomorphic
vector bundles on $\Si$ (see \cite{H} for references).

 The form $F_{A_{\CS}}$ (\ref{eq3.32}) is a $(1,1)$-form and
the line bundle $L$ admits a unique holomorphic structure compatible with
the Hermitian connection $A_{\CS}$. Moreover, a complex structure $I$ on
$\Si$ defines a K\"ahler Weil--Petersson metric on $\sM\semis$ with
K\"ahler form $\om_{\mathrm{WP}}=\Om$. This metric defines the
Levi-Civita connection on the complex tangent bundle $T \sM\semis$, and
hence a Hermitian connection $A_{\LC}$ on the line bundle
 \begin{equation}
\det T \sM\semis=L^{\tensor 4},
 \label{eq3.34}
 \end{equation}
and a Hermitian connection $A_{1/\LC}$ on $L$ compatible with the
holomorphic structure on $L$. Thus $A_{\LC}=A_{4\CS}$ and the equality
(1.3) holds with $d=-2$. We can use the geodesic lifting (2.23).

The result of complex quantization of the prequantum system (\ref{eq3.33})
can be viewed as the space of wave functions of level $k$, that is, the
space of $I$-holomorphic sections
 \begin{equation}
 \sH_I^k=H^0 (R_I, L^{k-2})
 \label{eq3.35}
 \end{equation}

One knows that this system of spaces and monomorphisms is related to the
system of representations of $\fsl(2,\C)$ in the
Weiss--Zumino--Novikov--Witten model of CQFT. The ranks of these spaces
are given by the Verlinde formula (see \cite{B}):
 \begin{equation}
 \rk \sH_I^k=\frac{k^{g-1}}{2^{g-1}} \sum_{n=1}^{k-1}
 \frac{1}{(\sin(\frac{n\pi}{k}))^{2g-2}}\,
 \label{eq3.36}
 \end{equation}
(please note the shift $k\mapsto k-2$). Many beautiful features of the
geometry of embeddings (1.7)
 \begin{equation}
 \PP \fie_k\colon R_I=\sM\semis\to \PP H^0(L^k)^*,
 \label{eq3.37}
 \end{equation}
 observed by Beauville, Laszlo, Pauly, Sorger and many others, make it
reasonable to call this area of mathematics the theory of {\em non-Abelian
theta functions}. But we would like to mention specially the observation
of Oxbury and Ramanan about the spaces of level 4 \cite{O}. In this case
the space (3.35) is the natural direct sum of spaces of Abelian theta
functions of Jacobian and Pryms of a Riemann surface and the union of
classical theta bases (3.3) is the Bohr--Sommerfeld basis for this case.

 Moreover the vector bundle (1.2)
 \begin{equation}
 \sH^k\to M_g
 \label{eq3.38}
 \end{equation}
 over the moduli space of Riemann surfaces of genus $g$ admits the
projectively flat Hitchin connection \cite{H}. 

Our space $R_g$, with complex structure induced by a complex structure $I$
on $\Si$, is a singular algebraic variety $R_I$ and
 \begin{equation} 
 \Sing R_I=K_\Si
 \label{eq3.39}
 \end{equation}
is the Kummer variety of $\Si$. The restriction 
 \begin{equation}
L^k\rest{K_\Si}=L_{\Theta}^{2k}.
\label{eq3.40}
 \end{equation}
It is easy to see that if $k\gg0$, the restriction gives the epimorphism
 \begin{equation}
 \res\colon H^0(L^k)\to H^0(L^{2k}_\Theta)\to 0
\label{eq3.41}
 \end{equation}
and using our Hermitian structure (1.8) we get the orthogonal decomposition
 \begin{equation}
 H^0(L^k)=H^0(L^{2k}_\Theta) \oplus H^0 (J_{K_\Si}\tensor L^k). 
\label{eq3.42}
 \end{equation}
For the first component of this decomposition we still have the special
basis (3.29), (3.30). Moreover this component decomposes as
 \begin{equation}
 H^0(L^{2k}_\Theta)=H^0(\Oh_{(J_\Si)_2})\oplus H^0 (J_{(J_\Si)_2}\tensor
L^k). 
\label{eq3.43}
 \end{equation}
Summarizing, we have a filtration of the vector bundle (3.38):
 \begin{equation}
 \sH_{\triv}\subset \sH^{2k}_{\red}\subset \sH^{k+2}
 \label{eq3.44}
 \end{equation}
corresponding to the flag (3.31) and the decompositions (3.42) and (3.43)
(see also (3.27) and (3.30)). Every bundle of this flag has a projective
flat connection, these connections are hereditary. The monodromies of the
first pair of bundles 
 \begin{equation}
 \Mon \sH_{\triv}=\Sp (2g, \Z_2) \quad\text{and}\quad
 \Mon \sH^{2k}_{\red}=\Sp (2g,\Z_{2k})
 \label{eq3.45}
 \end{equation}
are finite. Our main result (see 4.18) reduces the question of the
monodromy of $\sH^{k+2}$ to purely combinatorial question about the
representation (4.15). This representation is the subject of an absolutely
different and very beautiful theory providing many 3-manifold invariants.

\section{Combinatorial theory and identifications}

 Let $\Ga$ be any 3-valent graph having vertices $V(\Ga)$ and edges
$E(\Ga)$, with $|V(\Ga)|=2g-2$ and $|E(\Ga)|=3g-3$; consider functions
 \begin{equation}
 w\colon E(\Ga)\to \left\{0, \frac{1}{2k},\dots,\frac{1}{2}\right\}
 \label{eq4.1}
 \end{equation}
on the edges of $\Ga$ to rational numbers with denominator $\frac{1}{2k}$
in $[0,\half]$ satisfying:
 \begin{enumerate}
 \addtocounter{enumi}{-1}
 \item $w(C_i)\in \Z\cdot\frac{1}{k}$ if $C_i$ disconnects $\Ga$;
 \item[] and for any three edges $C_l,C_m,C_n$ meeting at a vertex $P_i$:
 \item $w(C_l) + w(C_m) + w(C_n)\in \frac{1}{k}\cdot \Z$;

 \item $w(C_l) + w(C_m) + w(C_n) \le 1$;

 \item for any ordering of $C_l, C_m, C_n$,
 \begin{equation}
 |w(C_l)-w(C_m)|\le w(C_n)\le w(C_l)+w(C_m).
 \label{eq4.2}
 \end{equation}

 \end{enumerate}
A function $w$ satisfying these conditions is called an {\em admissible
integer weight} of level $k$ on $\Ga$. Let $W^k_g(\Ga)$ be the set of
admissible integer weights. This set is canonically embedded in the set
$W_g^k$ (4.1) of all (unrestricted) functions 
 \begin{equation}
 W^k_g(\Ga)\subset W^k_g,
 \label{eq4.3}
 \end{equation}
which are obviously $|W^k_g|=(k+1)^{3g-3}$ in number.

 \begin{prop}[see for example \cite{K1}]
 The number\/ $|W_g^k(\Ga)|$ of admissible weights of level\/ $k$ is
independent of\/ $\Ga$.
 \end{prop}

The restrictions (\ref{eq4.2}) are called the {\em Clebsch--Gordan
conditions}. We can consider the space of all real functions with values in
$[0,1]$ subject to these conditions to get a complex $\De_{\Ga}$ (see
\cite{JW1}). We thus have a space
 \begin{equation}
\sH_{\Ga}^k=\bigoplus_{w\in W_g^k(\Ga)} \C \cdot w.
 \label{eq4.4}
 \end{equation}
with a natural Hermitian pairing $\Span{\ \,,\ }_c$ such that the $\{w\}$
form a unitary orthonormal basis. All these spaces are of course
canonically contained in the common space
 \begin{equation}
\sH_{\Ga}^k\subset \sH_g^k=\bigoplus_{w\in W_g^k} \C \cdot w.
\label{eq4.5}
 \end{equation}
The geometry of the projective configuration
 \[
 \bigcup_{\text{all graphs}} \PP \sH_{\Ga}^k\subset \PP \sH_g^k
 \]
reflects properties of the Moore--Seiberg complex (\cite{MS}) for 3-valent
graphs.

This combinatorial description has the following geometric meaning:
consider $\R^{3g-3}$ with coordinates $c_i$ corresponding to
$\{C_i\}=E(\Ga)$. It contains the complex $\De_{\Ga}$ and the integer
sublattice $\Z^{3g-3}\subset \R^{3g-3}$, and we can consider the
``action'' torus:
 \begin{equation} T^A=\R^{3g-3} / \Z^{3g-3}.
\label{eq4.6}
 \end{equation}
Then $T^A$ contains a topological complex $\Debar_{\Ga}$ obtained by
glueing together the boundary points of the polytope $\De_{\Ga}$.

Of course every unrestricted function $w\in W_g^k$ (4.1) defines a
$2k$-torsion point $w\in T^A_{2k}$ on the action torus. Thus we have an
identification:
 \begin{equation}
 W_g^k=T^A_{k}
 \label{eq4.7}
 \end{equation}

In particular, the admissible integer weights $W_g^k$ can be viewed as a
subset of the $2k$-torsion points of the action torus:
 \begin{equation}
 W_g^k(\Ga)\subset T^A_{k}.
 \label{eq4.8}
 \end{equation}

Now if we pump up the edges of a trivalent graph $\Ga$ to tubes, and the
vertices to small 2-spheres we get a Riemann surface $\Si_{\Ga}$ of genus
$g$ marked with $3g-3$ disjoint, noncontractible, pairwise nonisotopic
smooth circles $\{C_i\}$ on $\Si$, the meridian circles of the tubes. The
isotopy class of such a collection of circles is called a {\em marking} of
the Riemann surface. It is easy to see that the complement is the union
 \begin{equation}
 \Si_g-\{C_1, \dots, C_{3g-3}\}=\coprod_{i=1}^{2g-2} P_i
 \label{eq4.9}
 \end{equation}
of $2g-2$ trinions $P_i$, where every trinion is a 2-sphere with 3
disjoint discs deleted:
 \begin{equation}
 P_i=S^2 \setminus \bigl(D_1\cup D_2\cup D_3\bigr) \quad \text{with}
 \quad\Dbar_i\cap\Dbar_j=\emptyset \quad \text{for} \quad i \ne j.
 \label{eq4.10}
 \end{equation}

On the other hand any trinion decomposition of our Riemann surface $\Si$,
given by a choice of a maximal collection of disjoint, noncontractible,
pairwise nonisotopic smooth circles on $\Si$. It is easy to see that any
such system contains $3g-3$ simple closed circles
 \begin{equation}
C_1, \dots, C_{3g-3}\subset \Si_g,
 \label{eq4.11}
 \end{equation}
with complement the union of $2g-2$ trinions $P_j$. The type of such a
decomposition is given by its {\em $3$-valent dual graph} $\Ga(\{C_i\})$,
associating a vertex to each trinion $P_i$, and an edge linking $P_i$ and
$P_j$ to a circle $C_l$ such that
 \[
 C_l\,\subset\,\partial P_i\cap\partial P_j.
 \]
Thus the isotopy class of a trinion decomposition is given by a 3-valent
graph $\Ga$.

Now the modular group $\Mod_g$ which acts on $R_g$ by symplectomorphisms
preserving the prequantum data. Every element $\ga\in\Mod_g$ changes the
system of loops $\{[C_i]\}\to\{\ga([C_i])\}$ but the graph of the trinion
decomposition is precisely the same:
 \begin{equation}
 \Ga(\{[C_i]\})=\Ga(\{\ga([C_i])\}).
\label{eq4.12}
 \end{equation}
Thus the set of admissible integer weights $W_g^k(\Ga(\{[C_i]\}))$ is the
same, and defines the basis (4.4) in the space $\sH^k_{\Ga}$.

Moreover, using the fusion matrices that describe the monodromy of the
Knizhnik--Zamolodchikov equation, Kohno \cite{K1} constructed a canonical
isomorphism
 \begin{equation}
 \sH^k_{\Ga_1}=\sH^k_{\Ga_2}
 \quad\text{ for any two graphs $\Ga_1$ and $\Ga_2$,}
\label{eq4.13}
 \end{equation}
 and as a consequence, he obtained unitary linear representations of the
central extensions
 \begin{equation}
 1\to Z(k)\to \wMod_g^k\to \Mod_g\to 1,
 \label{eq4.14}
 \end{equation}
where $Z(k)$ is the cyclic group generated by $\exp(2\pi
i\frac{k}{8(k+2)})$; Kohno's representations are
 \begin{equation}
 \rho_{c}^k\colon \wMod_g^k\to \End \sH_I^k.
 \label{eq4.15}
 \end{equation}

 The decomposition with these components is parallel to the decomposition
of the highest weight representation of the affine Lie algebra of
$\fsl(2,\C)$ by eigenspaces of the operator $L_0$ from Sugawara's
construction of the representation of the Virasoro Lie algebra (see
\cite{K2}).

 Using this representation, we construct the vector bundles
 \begin{equation}
\pi\colon \mathcal H^k_c\to M_g
 \label{eq4.16}
 \end{equation}
(the subscript c stands for ``combinatorial'') over the moduli space $M_g$
of curves of genus $g$ having fibres
 \begin{equation}
\pi\1(\Si_{\Ga})=\sH^k_{\Ga},
 \label{eq4.17}
 \end{equation}
with the projective unitary connection $a_{c}$ with the monodromy
representation (4.15). Indeed, $\Mod_g$ acts transitively both on all
markings $\{[C_i]\}$, and on all trinion decompositions.

\subsection*{The main result} We want to identify the projectivizations of
the spaces (3.35) and (4.4):

 \begin{thm} For $k\gg0$ (depending on the genus $g$) then there exists a
canonical identification
 \begin{equation}
\sH^{k+2}_I=\sH^k_{\Ga_0}
 \label{eq4.18}
 \end{equation}
up to finite ambiguity for the special\/ $3$-valent graph\/ $\Ga_0$
described below. \end{thm}

 This identification gives a chain of identifications of objects and
construction of two theories: the complex quantization of (3.33) (the WZNW
model of CQFT) and the combinatorial theory of Witten, Reshetikhin--Turaev,
Tsuchiya--Kanie, Drinfeld, Moore--Seiberg and Kohno. Thus we can use
results of the theory of non-Abelian theta functions in algebraic
geometric as an effective means of computing topological invariants of
3-manifolds. On the other hand, the standard bases in the spaces
$\sH^k_{\Ga}$ (4.4) define non-Abelian theta functions {\em with
characteristic} of level $k$ and, following Mumford, we should write down
special equations in these bases defining the images of $\sM\semis$ in
spaces of conformal blocks.

 We realize this program using the BPU method of Sections~1--2. (Thus we
{\em only} get the identification (4.18) for $k\gg0$.)

We must construct $\BS^{k+2}$ cycles on $R_g$ indexed by the set
$W^k_g(\Ga)$ (4.3). This was done by Jeffrey and Weitsman \cite{JW1}:
for a marked Riemann surface $\Si_{\Ga}$, the map
 \begin{equation}
 \pi_{\{C_i\}}\colon R_g\to \R^{3g-3}
 \label{eq4.19}
 \end{equation}
with fixed coordinates $(c_1,\dots, c_{3g-3})$ such that
 \begin{equation}
 c_i (\pi_{\{C_i\}} (\rho))=\frac{1}{\pi}\,
\cos\1\bigl(\textstyle{\frac{1}{2}}\tr \rho([C_i])\bigr)\in [0, 1],
 \label{eq4.20}
 \end{equation}
where $\{C_i\}=E(\Ga)$. It is well known that
 \begin{enumerate}
 \item The map $\pi_{\{C_i\}}$ is a real polarization of the system
$(R_g,k\cdot\om,L^k,k\cdot A_{\CS})$.

 \item The coordinates $c_i$ are action coordinates for this
Hamiltonian system.
 \item The map $\pi_{\{C_i\}}$ is a {\em moment map} for the action of
$T^{3g-3}$ on $R_g$
 \begin{equation}
 R_g \times T^{3g-3}\to R_g
 \label{eq4.21}
 \end{equation}
constructed by Goldman \cite{G}.

 \item The image of $R_g$ under $\pi_{\{C_i\}}$
 is a convex polyhedron 
 \begin{equation}
 \De_{\{C_i\}}\subset [0,1]^{3g-3}. 
 \label{eq4.22} \end{equation}

 \item The symplectic volume of $R_g$ equals the Euclidean volume of
$\De_{\{C_i\}} $:
 \[
 \int_{R_g}\om^{3g-3}=\Vol \De_{\{C_i\}}=\frac{2\cdot \zeta
(2g-2)}{(2\pi)^{g-1}}.
 \]

 \item The expected number of Bohr--Sommerfeld orbits of the real
polarization $\{C_i\}$
 \begin{equation}
 N_{\BS}(\pi_{\{C_i\}},R_g,\om,L,A_{\CS})
 \label{eq4.23}
 \end{equation}
is equal to the number of half-integer points in the polyhedron
$\De_{\{C_i\}}$, and
 \begin{equation}
 \lim_{k\to\infty} \frac{N_{\kBS}}{k^{3g-3}}\ =\ \int_{R_g}\om^{3g-3} =\Vol
\De_{\{C_i\}}.
 \label{eq4.24}
 \end{equation}
 \end{enumerate}
These functions $c_i$ are continuous on all $R_g$ and smooth over $(0,1)$.

Every $w\in W^k_g(\Ga)$ defines a point of $\De_{\{C_i\}}$ (4.22)
with coordinates
 \[
c_i=2w(C_i).
 \]

 \begin{prop}[\cite{JW1}] \begin{enumerate}
 \item The map $x\mapsto 2x$ sends the complex $\De_\Ga$ (4.7) to
$\De_{\{C_i\}}$ (4.22)
 \[
2\De_\Ga=\De_{\{C_i\}};
 \]

 \item this transformation sends $W^k_g(\Ga)$ to the set of\/ $\BS_k$
fibres of the real polarization $\pi_\Ga=\pi_{\{C_i\}}$ (4.19):
 \begin{equation}
2 W^k_g(\Ga)=\De_{\{C_i\}}\cap \BS^{k+2}(L);
\label{eq4.25}
 \end{equation}
 
 \item in particular
 \[
 |W_g^k(\Ga)|= \#(\De_{\{C_i\}}\cap \BS^{k+2}(L)).
 \]
 \end{enumerate}
 \end{prop}

In summary:
 \begin{enumerate}
 \item fixing the graph $\Ga(\{[C_i]\})$ of a trinion decomposition we
enumerate canonically the set of $k$-Bohr--Sommerfeld fibres of all
polarizations with the same graph as the set $W_g^k(\Ga(\{[C_i]\}))$ of
integer weights on this graph;
 \item fixing a collection of loops $\{[C_i]\}$ we get a finite set of
disjoint $k$-Bohr--Sommerfeld oriented cycles $\sL_w$ for $w\in
W_g^k(\Ga(\{[C_i]\}))$ in $R_g$;
 \item for any level $k$, any complex Riemann surface $\Si$, and any
trinion decomposition $\{C_i\}$, we have
 \begin{equation}
 \rk \sH_I^{k+2}=\rk \sH^k_{\Ga}=\rk\sH^k_{\pi}
 =\text{Verlinde number (3.36)}.
 \label{eq4.26}
 \end{equation}
 \end{enumerate}

We complete the description in \cite{JW1} of the $\BS^{k+2}$ fibres of a
real polarization $\pi_{\Ga}$ by describing the fibres $\pi\1_\Ga(w)$ for
which
 \begin{equation}
 \pi\1_\Ga(w)\cap K_\Si \ne\emptyset.
 \label{eq4.27}
 \end{equation}

 \begin{rmk} We will see below that $\BS^{k+2}$ fibres disjoint from
$K_\Si=\Sing R_I$ can only have orbifold singularities. So we can apply the
BPU construction to it, and get a partial basis of ``theta functions with
characteristic''.
 \end{rmk}

Return to the geometric procedure described after formula (4.8). Pumping up
our graph $\Ga$ we get a handlebody $\wGa$ with boundary
 \begin{equation}
 \partial\wGa=\Si,
 \label{eq4.28}
 \end{equation}
giving an exact sequence of fundamental groups
 \[
1\to\ker\to \pi_1(\Si)\to \pi_1 (\wGa)\to 1,
 \]
where the kernel is the subgroup of the fundamental group of Riemann
surface of cycles homotopic to a point in $\wGa$. To recognize our previous
handlebody $\wSi_a$ (3.17), recall that the standard presentation (3.16) of
the fundamental group of $\Si_g$ defines another ``dual'' presentation
given by the following GNW construction (Guruprasad--Nilakantan--Weil, see
\cite{T4} for references): set 
 \[
 \al_i=r_{i-1} b_i\1 r_i\1, \quad \beta_i=r_i a_i\1 r_{i-1}\1,
 \quad\text{where}\quad r_i=\prod_{j=1}^{i} [a_j, b_j].
 \]
Then
 \[
\pi_1(\Si_g)=\Span{\al_1,\dots, \al_g, \beta_1,\dots,\beta_g
\bigm|\prod\nolimits_{j=1}^{g} [\al_j, \beta_j]=1}
 \]
is another presentation of $\pi_1(\Si_g)$. Sending the generators $a_i,
b_j$ to $\al_i, \beta_j$ gives an automorphism $W$ of $\pi_1(\Si_g)$, that
is, $W\in\Mod_g$, and it is an involution: $W^2=\id$.

To show that
 \begin{equation}\wGa=\wW(\Si)_a=\wSi_{W(a)},
\label{eq4.29}
 \end{equation}
 consider the special 3-valent graph $\Ga_0$ corresponding to the
presentation $\Si$ as a connected sum of $g$ 2-tori (see \cite{K1},
Figs.~12a and~13b). We get a basis (3.16). To get $\Ga_0$, we fix the
following system of cycles $\{C_i\}$ on $\Si$: they consist of three
groups:
 \begin{enumerate}
 \item $a_1,\dots,a_g$, the cycles $a_i$ of \cite{K1}, Fig.~13b;
 \item $a'_2,\dots,a'_{g-1}$, the cycles $c_i$ of \cite{K1}, Fig.~13b;
 \item $c_1=[a_1,b_1],\dots,c_{g-1}=[a_{g-1},b_{g-1}]$, the cycles $b_i$
of \cite{K1}, Fig.~13b.
 \end{enumerate}
 Then 
 \begin{equation}
\{c_i\}\subset [\pi_g, \pi_g]
\label{eq4.30}
 \end{equation}
is the commutator subgroup of the fundamental group. Since $E(\Ga_0)$
contains the subset
 \begin{equation}
E(\Ga_0)^a=\{a_1,\dots, a_g\}.
\label{eq4.31}
 \end{equation}
Our handlebody (3.17) transformed by $W$ is equal to $\wGa_0$.

We label the coordinates $\{c_i\}$ (4.19) and (4.20) by the same symbols
$a_i,a'_j,c_k$. Then 
 \begin{align}
 \rho\in K_\Si &\implies c_i (\rho)=0 \quad \text{for} \quad i=1, 2,\dots,
 g-1; \label{eq4.32} \\
 \rho\in K_\Si &\implies a_i=a'_i \quad \text{for}\quad i=2,\dots, g-1.
\label{eq4.33} 
 \end{align}
Thus
 \[
\dim \pi_{\Ga_0}(K_\Si)=g;
 \]
more precisely, we have:
 \begin{prop}
 \[
 \pi_{\Ga_0}(K_\Si)=K_{\wGa_0};
 \]
and $\pi_{\Ga_0}$ is the real polarization (3.22).
 \end{prop}

We must now check the following:

 \begin{prop}
 \[
 K_{\wGa_0}\cap \BS^k(R_g, L)=K_{\wGa_0}\cap \BS^{2k}(K_\Si, L_\Theta)
 \]
 \end{prop}

The proof follows immediately from the description of $\BS_k$ fibres in
Proposition~5.2 and Corollary~5.1 of the next section.

 \begin{cor}
 \begin{enumerate}
 \item $W^k_g(\Ga_0)$ contains the subset
 \[
W^k_g(\Ga_0)\Ab=K_{\wGa_0}\cap \BS^{k+2}(L)\subset W^k_g(\Ga_0)
 \]
of weights that we call\/ {\em Abelian} weights.

 \item Weights of the set 
 \[
 W^k_g(\Ga_0)\nonAb=W^k_g(\Ga_0)\setminus W^k_g(\Ga_0)\Ab
 \]
are called\/ {\em non-Abelian weights}.

 \item The space (4.4) can be decomposed into Abelian and non-Abelian parts
 \[
\sH_{\Ga_0}^k=\Biggl(\bigoplus_{w\in W^k_g(\Ga_0)\Ab} \C \cdot w\Biggr)
\oplus \Biggl(\bigoplus_{w\in W^k_g(\Ga_0)\nonAb} \C \cdot w\Biggr).
 \]

 \item (3.30) identifies the Abelian component:
 \[
 \Biggl(\bigoplus_{w\in W^k_g(\Ga_0)\Ab} \C \cdot w\Biggr)
 =H^0(K_\Si, L^{2k}_\Theta).
 \]
 \end{enumerate}
 \end{cor}

To get a basis of ``theta functions with characteristic'' in all
$H^0(L^k)$ related to a trinion decomposition of $\Si$ and a linear
isomorphism
 \[
\BPU_k\colon \mathcal H_{\pi_{\Ga_0}}^k\to H^0(L^k)
 \]
(2.24) we must construct on every fibre $\sL_w$ for $w\in W_g^k$ an almost
canonical half-form $\hF_w$ in such a way that the Szeg\"o projector (2.2)
extends to a class of distributions including $(\La_w,\hF_w)$ for every
$w\in W_g^k(\Ga)$.

\section{Covariant constant half-forms and\\ singularities}

To put covariant constant half-forms on BS fibres, recall some facts about
its structure. Let $\Ga$ be a 3-valent graph with vertices $V(\Ga)$ and
edges $E(\Ga)$, and suppose that
 \[
w\colon E(\Ga)\to \left\{0, \frac{1}{2k},\dots, \frac{1}{2}\right\}
 \]
is an integer admissible weight. For $\al\in \{0, \frac{1}{2k},\dots,1\}$
let
 \begin{equation}
w\1(\al)=\Ga_1(\al)\cup\cdots\cup\Ga_n(\al)\subset\Ga
\label{eq5.1}
 \end{equation}
be the decomposition into connected components. Then every component
$\Ga_i(\al)$ is a 3-valent graph with $n_i$ univalent vertices $a_1,\dots,
a_{n_i}$ (see \cite{K1}).

Every gauge class of connections contains a connection $a_0$ adapted to a
trinion decomposition (see \cite{JW1}, Definition~2.2). Fix the
filtration
 \begin{equation}
Z(\SU(2))=\Z_2\subset \U(1)\subset \SU(2),
\label{eq5.2}
 \end{equation}
and view it as the triple
 \begin{equation}
G=\{\Z_2, \U(1), \SU(2)\}
\label{eq5.3}
 \end{equation}
For $[a]\in \pi\1(w)$, we have the function
 \begin{equation}
e_w\colon E(\Ga)\to G
\label{eq5.4}
 \end{equation}
 sending every loop $C_j$ to the element of $G$ conjugate to the
stabilizer of the monodromy of $[a]$ around this loop, and the function
 \begin{equation}
v_w\colon V(\Ga)\to G
\label{eq5.5}
 \end{equation}
sending a trinion $P_i$ to the stabilizer of the flat connection
$a\rest{P_n}$. Of course,
 \begin{align}
 & C_j\subset \partial P_n \implies v_w(P_n)\subset e_w(C_j); \notag \\
 & C_1 \cup C_2 \cup C_3=\partial P_n \quad \text{and} \quad
e_w(C_1)=e_w(C_2)=\SU(2) \implies 
 \label{eq5.6}
 \\
 & \qquad\qquad e_w(C_3)=\SU(2) \implies v_w(P_n)=\SU(2), \notag
 \end{align}
and so on.

Obviously
 \begin{equation}
e_w(C_j)= \text{$\U(1)$ or $\SU(2)$.}
\label{eq5.7}
 \end{equation}

 \begin{prop} The functions $e_w$ and $v_w$ depend on only $w$ and not on
the choice of $[a]\in \pi\1(w)$.
 \end{prop}

More precisely, they depend on the combinatorics of the decomposition
(5.1).

Thus $w$ defines direct products
 \[
 \prod_{C\in E(\Ga)} e_w(C) \quad\text{and}\quad
 \prod_{P\in V(\Ga)} v_w(P),
 \]
and $\prod_{P\in V(\Ga)} v_w(P)$ acts on $\prod_{C\in E(\Ga)} e_w(C)$ as
follows: for
 \begin{gather*}
 g=(g_1,\dots, g_{2g-2})\in \prod_{P\in V(\Ga)} v_w(P)
 \quad\text{with $g_i\in v_w(P_i)$, and}\\
(t_1,\dots, t_{3g-3})\in \prod_{C\in E(\Ga)} e_w(C)
 \quad\text{with $t_n\in e_w(C_n)$,}
 \end{gather*}
 if $C_n\subset \partial P_i\cap\partial
P_j$ then
 \begin{equation}
 g(t_n)=g_i \circ t_n \circ g_j\1.
 \label{eq5.8}
 \end{equation}

 \begin{prop}[\cite{JW1}, Theorem~2.5] The fibre $\pi\1(w)$ is given by
 \begin{equation}
 \pi\1(w)=\prod_{C\in E(\Ga)} e_w(C) \Bigm/ \prod_{P\in V(\Ga)} v_w(P).
 \label{eq5.9}
 \end{equation}
 \end{prop}

Applying this description to
 \[
w\in K_{\wGa_0}\cap \BS^k(R_g, L)
 \]
proves Proposition 4.4: in this case, for every $P\in V(\Ga_0)$
 \[
 \U(1)\subset v_w(P) \quad \text{and} \quad e_w(c_i)=\SU(2)
 \]
for $c_i$ from (4.30). 

 \begin{cor}
 The fibre $\pi\1(w)$ is isomorphic to
 \begin{equation}
 \pi\1(w)\iso T^t \times [(S^3)^p \times (S^2)^s] / G_w,
 \label{eq5.10}
 \end{equation}
where $t, p$ and $s$ are nonnegative integers and $G_w$ is the finite
Abelian group defined by $w$, or more precisely by the combinatoric data
(4.1); moreover,
 \begin{equation}
 H_1(\pi\1(w))=\Z^t \oplus \Z_2^p.
 \label{eq5.11}
 \end{equation}
 \end{cor}

Translations along the torus $T^t$ in (4.10) are induced by Hamiltonians
lifted from the target space $\R^{3g-3}$ of $\pi$. We consider below
half-forms invariant under such translations.

Jeffrey and Weitsman \cite{JW1} use the normalization of the action
coordinates via branched covers to construct a covariant constant section
$s_w$ of the restrictions of $(L^k,A_{k\text{-CS}})$ to $\pi\1(w)$.

Our groups (5.2--3) admit bi-invariant half-forms $\hF_1$ on $\U(1)$ and
$\hF_3$ on $\SU(2)$. For every $w$ we can normalize these form $\hF_1(w)$
and $\hF_2(w)$ so that the half-form
 \begin{equation}
 \hF_w=(\hF_1(w))^{t-s} \cdot (\hF_2(w))^{p+s}
 \label{eq5.12}
 \end{equation}
is homogeneous of degree 1 on $\pi\1(w)$ (see (5.10)) with respect to
scaling $\hF_i(w)\to t\cdot\hF_i(w)$. We say that such half-form is
homogeneous normalized.

It's easy to see (\cite{JW2}, 4.7) that a normalized half-forms for a
nonsingular $\BS^{k+2}$ fibre is given by a Hamiltonian vector field with
Hamiltonian in $\R^{3g-3}$ of volume 1.

Thus every $\BS^{k+2}$ fibre is given the covariant constant half-form
(5.12), and we can proceed to construct the corresponding Legendrian
distributions in $P$. Recall that $R_g$ is almost homogeneous with respect
to the Goldman torus action (4.21). Thus the Schwartz kernel of coherent
states does not depend on points and, outside singular points of fibres,
they behave as in the homogeneous case (see \cite{BPU}, (10--13)). By
lifting to $P$ every $\BS^{k+2}$ fibre $\pi\1(w)$ defines Legendrian
subcycle $\La_w\subset P$ marked with the half-form
 \begin{equation}
 \hFbar_w=(\fie_4\circ {\det}\circ G(i))^*\hF_w,
 \label{eq5.13}
 \end{equation}
and having monodromy a $k$th root of 1.

To apply the BPU construction the principal bundle $P=S^1(L^*)$ must be
given a metaplectic structure. This can be done at once using (3.34).

For the identification (4.18), consider the decompositions (3.42)
 \[
 H^0(L^k)=H^0(K_\Si, L^{2k}_\Theta) \oplus H^0(J_{K_\Si}\otimes L^k)
 \]
and 
 \[
 \sH_{\Ga_0}^k=\Biggl(\bigoplus_{w\in W^k_g(\Ga_0)\Ab} \C \cdot w\Biggr)
 \oplus \Biggl(\bigoplus_{w\in W^k_g(\Ga_0)\nonAb} \C \cdot w\Biggr) 
 \]
of Corollary 4.1. Then (3.30) identifies the first (Abelian) components and
the BPU construction identifies the second (non-Abelian) components. This
method (and its verbatim modification for orbifolds) is applicable
because the non-Abelian $\BS_k$ fibres are contained in the smooth part of
$R_g$.

\subsection*{Acknowledgments}

I would like to express my gratitude to the Mathematics Institute of
Warwick University and personally to Miles Reid and Victor Pidstrigatch for
support and hospitality. I was very influenced by the work of John
Rawnsley. Let him be thanked for this. Thanks are due to Miles Reid for
efforts to make this paper readable.

\bigskip
\noindent
Andrei Tyurin, Algebra Section, Steklov Math Institute,\\
Ul.\ Gubkina 8, Moscow, GSP--1, 117966, Russia \\
e-mail: Tyurin@tyurin.mian.su {\em or} Tyurin@Maths.Warwick.Ac.UK\\
{\em or} Tyurin@mpim-bonn.mpg.de
 \end{document}